\documentclass{amsart}

\usepackage{amssymb,amscd,amsthm,amsxtra}
\usepackage{dsfont,latexsym}

\newtheorem{thm}{Theorem}[section]
\newtheorem{corollary}[thm]{Corollary}
\newtheorem{lemma}[thm]{Lemma}
\newtheorem{prop}[thm]{Proposition}
\theoremstyle{definition}

\theoremstyle{remark}
\newtheorem{rem}[thm]{Remark}
\numberwithin{equation}{section}

\newcommand{\eps}{\varepsilon}

\newcommand{\st}{\,{\mbox{ s.t. }}\,}

\newcommand{\R}{\mathds{R}}
\newcommand{\N}{\mathds{N}}
\newcommand{\Z}{\mathds{Z}}

\newcommand{\CC}{\mathcal{C}}
\newcommand{\FF}{\mathcal{F}}
\newcommand{\EE}{\mathcal{E}}
\newcommand{\KK}{\mathcal{K}}

\usepackage{mathrsfs}
\renewcommand{\mathcal}{\mathscr}

\renewcommand{\le}         {\leqslant}

\renewcommand{\ge}         {\geqslant}

\renewcommand{\emptyset}         {\varnothing}

\begin{document}

\author{Ovidiu Savin and Enrico Valdinoci}
\thanks{OS has been supported by 
NSG grant 0701037. EV has been
supported
by
FIRB 
project ``Analysis and Beyond''
and GNAMPA
project
``Equa\-zio\-ni non\-li\-nea\-ri su va\-rie\-t\`a:
pro\-prie\-t\`a qua\-li\-ta\-tive e clas\-si\-fi\-ca\-zio\-ne
del\-le so\-lu\-zio\-ni''.
Part of this work
was carried out
while EV was visiting Columbia
University.}

\title[A variational model
driven by the
Gagliardo norm]{Density estimates for a variational model \\
driven by the Gagliardo norm}

\begin{abstract}
We prove density estimates for level sets of minimizers of the energy 
$$\eps^{2s}\|u\|_{H^s(\Omega)}^2+\int_\Omega W(u)\,dx,$$ with $s \in 
(0,1)$, 
where $\|u\|_{H^s(\Omega)}$ denotes the total contribution from $\Omega$ 
in the $H^s$ norm of $u$, and $W$ is a double-well potential. 

As a consequence we obtain, as $\eps \to 0^+$,
the uniform convergence of the level sets of $u$ to either a 
$H^s$-nonlocal minimal surface if $s\in(0,\frac 1 2)$, or to a classical 
minimal surface if $s \in[\frac 1 2,1)$.
\end{abstract}

\maketitle


\section{Introduction}

A classical model for the energy of a two-phase fluid of density $u$ lying 
in a bounded
domain $\Omega \subset \mathbb{R}^n$, with 
$n\ge2$,
is given by the Ginzburg-Landau energy functional $$\int_\Omega \frac{\eps^2}{2}|\nabla u|^2+W(u) \, dx.$$ The function $W:\R\rightarrow [0,+\infty)$ is a double well potential with two zeros (minima) at the densities of the stable phases, which we assume for simplicity to be $+1$ and $-1$. The kinetic energy is given by the Dirichlet integral $$ \frac{\eps^2}{2} \int |\nabla u|^2 dx$$ which takes into account interactions at small scales between the fluid particles. The typical energy minimizer has two regions where $u$ is close to $+1$ and $-1$ which are separated by a ``phase transition" which lies in an $\eps$ neighborhood of the $0$ level set $\{u=0\}$.

In this paper we consider a different model in which the kinetic term is 
replaced by the $H^s$ (semi)norm of $u$, i.e.
$$\eps^{2s} \|u\|_{H^s}^2 
\quad \mbox{ 
with $s \in (0,1)$}.$$ This means that the interactions at small scales 
have nonlocal character. In this case the boundary data for $u$ is defined 
in $\mathcal {C} \Omega$, that is the complement of~$\Omega$. Similar 
models driven by
a fractional, Gagliardo-type norm were considered in \cite{Gag1,Gag2}; see 
also~\cite{ABS,GP,Gon}
and references therein
for a onedimensional related system
that models phase transitions on an interval.

{F}rom the physical point of view, the importance of
these type of models relies in their attempt to capture,
via the nonlocal term, the features arising from the long-range
particle interactions, and it is of course desirable to
understand if and how the
nonlocal aspect influences the interfaces
and to have good estimates on their width. Our results
are a first attempt to give some answers to
this questions. Indeed,
we show that the level sets of the minimizers for this nonlocal energy 
satisfy a uniform density property. For the Allen-Cahn-Ginzburg-Landau
energy such density estimates
were proved in~\cite{CC}.
As a consequence we obtain that when $s\in(0, 1/2)$ the phase 
transition 
converges locally uniformly as $\eps \to 0^+$ to a $H^s$-nonlocal minimal 
surface (see \cite{CRS} for the precise definition), and when $s
\in[1/2,1)$ the phase transition converges locally uniformly to a 
classical 
minimal surface.

We define
$$ X:=\big\{ u\in L^\infty (\R^n) \st \| u\|_{L^\infty(\R^n)}\le 
1
\big\},$$
the space of admissible functions -- when dealing with a minimization
problem in $\Omega$, we 
prescribe $u\in X$ with boundary data $u_o$ outside 
$\Omega$ (i.e., $u=u_o$ in $\CC\Omega$), and we say that a 
sequence~$u_n\in X$ converges to~$u$ in~$X$
if~$u_n$ converges to~$u$ in~$L^1(\Omega)$.

We define also
$$ \KK(u;\Omega):=\frac{1}{2}\int_\Omega\int_\Omega\frac{
|u(x)-u(y)|^2
}{
|x-y|^{n+2s}}\,dx\,dy
+\int_\Omega\int_{\CC\Omega}\frac{
|u(x)-u(y)|^2
}{
|x-y|^{n+2s}}\,dx\,dy,$$
the $\Omega$ contribution in the $H^s$ norm of $u$ $$\int_{\R^n} 
\int_{\R^n} \frac{
|u(x)-u(y)|^2
}{
|x-y|^{n+2s}}\,dx\,dy,$$ i.e we omit the set where $(x,y) \in \CC\Omega \times \CC \Omega$ since all $u\in X$ are fixed outside $\Omega$.

The energy functional $J_\eps$ in $\Omega$ is defined as  $$J_\eps(u;\Omega):=\eps^{2s}\KK(u;\Omega) + \int_\Omega W(u) \, dx.$$ Throughout the paper we assume that $W:[-1,1] \to [0, \infty)$,
\begin{equation}\label{Wcond}
 W \in C^2([-1,1]), \quad W(\pm 1)=0, \quad W>0 \quad \mbox{in $(-1,1)$}
 \end{equation}
$$W'(\pm 1)=0, \quad  {\mbox{and}}\quad
W''(\pm 1)>0.$$
We say that $u$ is a minimizer\footnote{A
similar notion of minimizer will hold, later on,
for the suitably rescaled versions of $J_\eps$, namely $\FF_\eps$
and $\EE$.}
for $J_\eps$ in $\Omega$ if 
$$J_\eps(u;\Omega) \le J_\eps(v;\Omega)$$ for any $v$ which coincides with 
$u$ in $\CC \Omega$.

We remark is that if $u$ minimizes $J_\eps$ in $\Omega$ then it minimizes $J_\eps$ in any subdomain $\Omega'\subset \Omega$ since 
\begin{eqnarray*}
&& \KK(u; \Omega)=\KK(u; \Omega')\\
&&+\frac{1}{2}\int_{\Omega-\Omega'}\int_{\Omega-\Omega'}\frac{
|u(x)-u(y)|^2
}{
|x-y|^{n+2s}}\,dx\,dy
+\int_{\Omega-\Omega'}\int_{\CC\Omega}\frac{
|u(x)-u(y)|^2
}{
|x-y|^{n+2s}}\,dx\,dy,\end{eqnarray*}
and the latter two
integral terms do not depend on the values of $u$ in 
$\Omega'$.

If~$u$ is a minimizer for~$J_\eps$ in all bounded open sets~$\Omega$,
we say, simply, that~$u$ is a minimizer\footnote{Sometimes,
in the literature, minimizers are called ``local'', or ``class A'',
minimizers.}.
The behavior of $J_\eps$ as $\eps \to 0^+$ is quite different as $s\in 
(0,1/2)$, $s=1/2$ and $s\in (1/2,1)$. We showed in \cite{SaV} that in each case $J_\eps$ must be multiplied by an appropriate constant depending on $\eps$ in order to obtain the $\Gamma$-convergence to a limiting functional. More precisely, given any~$\eps>0$, we define the
functional~$\FF_\eps:X\rightarrow\R\cup\{+\infty\}$ as
\begin{equation}\label{defu}\begin{split}&
\FF_\eps(u)=\FF_\eps(u;\Omega) :=\left\{
\begin{matrix}
\eps^{-2s} J_\eps(u; \Omega) & {\mbox{ if $s\in(0,\,1/2)$,}} \\
|\eps \log \eps|^{-1} J_\eps(u; \Omega) & {\mbox{ if $s=1/2$,}}\\
\eps^{-1} J_\eps(u; \Omega) & {\mbox{ if $s\in(1/2,\,1)$.}}
\end{matrix}
\right.\end{split}\end{equation}

In the case when $s\in (0,1/2)$, the limiting  functional ~$\FF:X\rightarrow\R\cup\{+\infty\}$
is defined as
\begin{equation}\label{FF1}
\FF(u):=\left\{\begin{matrix}
\KK(u;\Omega)
& {\mbox{ if $u|_\Omega = \chi_E -\chi_{\CC E}$, for some set $E\subset \Omega$}}
\\+\infty & {\mbox{otherwise.}}
\end{matrix}\right.
\end{equation}
In this case,~$\FF$
agrees with the nonlocal area functional of $\partial E$ in $\Omega$ that 
was studied in~\cite{CRS, CV, ADP}.
Remarkably, such nonlocal area functional is well defined
exactly when~$s\in(0,\,1/2)$.

In the case when $s \in [1/2,1)$ the limiting functional ~$\FF:X\rightarrow\R\cup\{+\infty\}$
is defined as
\begin{equation}\label{FF2}
\FF(u):=\left\{\begin{matrix}
c_\star \,{\rm Per}\,(E;\Omega)
& {\mbox{ if $u|_\Omega = \chi_E -\chi_{\CC E}$, for some set $E\subset \Omega$}}
\\+\infty & {\mbox{otherwise.}}
\end{matrix}\right.
\end{equation}
where $c_\star$ is a constant depending on $n$, $s$ and $W$.

We recall the $\Gamma$-convergence results in \cite{SaV}:

\begin{thm}\label{TH1}
Let~$s\in(0,\,1)$ and~$\Omega$ be a Lipschitz domain.
Then,
$\FF_\eps$ $\Gamma$-converges to~$\FF$, i.e.,
for any~$u\in X$,
\begin{itemize}
\item[{(i)}] for any~$u_\eps$ converging to~$u$ in~$X$,
$$ \FF(u)\le\liminf_{\eps\rightarrow0^+}\FF_\eps(u_\eps),$$
\item[{(ii)}] there exists~$u_\eps$ converging to~$u$
in~$X$ such that
$$ \FF(u)\ge\limsup_{\eps\rightarrow0^+}\FF_\eps(u_\eps).$$
\end{itemize}

\end{thm}

\begin{thm}\label{comp}
If $\FF_\eps(u_\eps; \Omega)$ is uniformly bounded for a sequence of $\eps 
\to 0^+$, then there exists a convergent subsequence $$u_\eps \rightarrow 
u_*:=\chi_E-\chi_{\CC E} \quad \mbox{ in $L^1(\Omega)$.}$$
Moreover, if $u_\eps$ minimizes $\FF_\eps$ in $\Omega$,
\begin{itemize}
\item[{(i)}] if~$s\in(0,\frac 1 2)$ and $u_\eps$ converges weakly to $u_o$ 
in $\CC 
\Omega$, then $u_*$ minimizes $\FF$ in~\eqref{FF1} 
among all the functions that coincide with~$u_o$ in~$\CC\Omega$;
\item[{(ii)}] if~$s\in[1/2,1)$, then $u_*$ minimizes~$\FF$ in 
\eqref{FF2}.
\end{itemize}
\end{thm}

Theorem~\ref{TH1} may be seen as
a nonlocal analogue of the celebrated $\Gamma$-convergence
result of~\cite{MM} (see also~\cite{Mo,Bou,OS} for further extensions). In
this framework, we recall
that a very important issue, besides $\Gamma$-convergence,
is the ``geometric'' convergence of the level sets of
minimizers to the limit surface. This topic has
been widely studied in the case of local functionals by
using appropriate density estimates (see~\cite{CC}, and
also~\cite{FV} and references therein for several other
applications). The idea of these density estimates is to give
an optimal bound on the measure occupied by the level sets
of a minimizer in a ball.

Our results give a nonlocal counterpart of these
density estimates for minimizers of $J_\eps$ (or $\FF_\eps$). For this, it is convenient to
scale space by a factor of $\eps^{-1}$ so that the dependence of $J_\eps$ 
on $\eps$ disappears. To be more precise, if $u$ minimizes $J_\eps$ in 
$\Omega$, then the rescaled function 
$$u_\eps(x):=u( \eps x)$$ minimizes $\EE $
in $\Omega_\eps:=\Omega / \eps$, where $$\EE(v; \tilde \Omega):= J_1(v;
\tilde \Omega)=\KK(v; \tilde \Omega) + \int_{\tilde \Omega} W(v) \, dx.$$

Our first result gives a uniform bound for the energy $\EE$ of a minimizer in $B_R$ for large $R$.

\begin{thm}\label{EB}
Let~$u$ be a minimizer of $\EE$ in~$B_{R+2}$ with $R \ge 1$. Then 
\begin{equation}\label{hole}
\EE(u;B_R)\, \le \,\left\{\begin{matrix}
\overline{C} \,R^{n-2s} & {\mbox{ if $s\in(0,\,1/2)$,}}\\
\overline{C} \,R^{n-1}\,\log R &{\mbox{ if $s=1/2$,}}\\
\overline{C} \,R^{n-1} & {\mbox{ if $s\in(1/2,\,1)$,}}
\end{matrix}\right.\end{equation}
where~$\overline{C}$ is a positive constant depending on~$n$,~$s$, and $W$.
\end{thm}

Theorem \ref{EB} can be stated in terms of minimizers~$u_\eps$ of 
$\FF_\eps$ in $B_{1+2\eps}$ as
$$\FF_\eps(u_\eps;B_1) \le \overline{C}.$$

Then, we have the following
density estimate on the level sets of minimizers:

\begin{thm}\label{DE}
Let~$u$ be a minimizer of $\EE$ in~$B_R$. Then for any~$\theta_1,$ $\theta_2\in(-1,1)$ such that
\begin{equation}\label{MUCO1}
u(0)>\theta_1, \end{equation}
we have that
\begin{equation}\label{hard}
\big| \{ u>\theta_2\}\cap B_R\big|\,\ge\, \overline{c} \,R^n
\end{equation}
if~$R\ge \overline{R}(\theta_1,\theta_2)$. The constant~$\overline{c}>0$ depends only on $n$, $s$ and $W$ and~$\overline{R}(\theta_1,\theta_2)$ is a large constant that depends also on $\theta_1$ and $\theta_2$.
\end{thm}

By scaling, Theorem \ref{DE} gives the
uniform density estimate for minimizers~$u_\eps$ of $\FF_\eps$ (or 
$J_\eps$) in 
$B_r$: if $u_\eps(0)> \theta_1$ then $$| \{ u_\eps>\theta_2\}\cap 
B_r|\,\ge\, 
\overline{c} \,r^n \quad \quad \mbox{if $r \ge \overline{R} \eps$.}$$

\begin{rem} {\rm
Our assumptions on $W$ are not the most general. For 
example in the 
Theorem \ref{EB} it suffices to say that~$W$ is bounded and $W(\pm 1)=0$. 
Also in Theorem \ref{DE} it suffices to assume that
there exists a small constant~$c>0$ such that
\begin{equation}\label{grow}
\begin{split}
&{\mbox{$W(t)\ge W(r)+c(1+r)(t-r)+c\, (t-r)^2$ when
$-1\le r\le
t\le
-1+c$}}\\
&{\mbox{and $W(r)-W(t)\le (1+r)/c$ when $-1\le r\le t\le+1$.}}
\end{split}
\end{equation}
Of course,~\eqref{grow} is warranted by our assumptions 
in~\eqref{Wcond},
but we would like to stress that less smooth, or even discontinuous,
potentials, may be dealt with using~\eqref{grow}.}\end{rem}

The proof of Theorem~\ref{DE}
is
contained in Section~\ref{SS2} and it requires
a careful analysis
of the measure theoretic properties of the minimizers
and several nontrivial modifications of the original proof
of~\cite{CC}, together with some iteration techniques
of~\cite{CRS}. In particular, the construction
of a new barrier function is needed in order to keep
track of the densities of the level sets in larger and
larger balls. 
Also, the proof
of~\eqref{hard} is somewhat delicate and it requires\footnote{When
$s\in(0,1/2)$ there is also an alternative
approach based on the fractional Sobolev inequality.
We will perform this different proof in~\cite{SaV-f}.}
the following
estimate for the double integral
\begin{equation}\label{1.8bis}
L(A,D):=\int_A\int_D 
\frac{1}{|x-y|^{n+2s}}dxdy.\end{equation}

\begin{thm}\label{GMT}
Let $s\in (0,1)$.
Let $A$ and $B$ be disjoint measurable subsets of $\R^n$
and let~$D:=\CC(A\cup B)$.
Then, there exists~$c\in(0,1)$, possibly depending on~$n$
and~$s$, for which the following
estimates hold:
\begin{itemize}
\item if $|B| \le c \,|A|$ and~$|A|>0$, then
\begin{equation}\label{gmt1}
L(A,D) \ge \left\{
\begin{matrix}
c\,|A|^{(n-2s)/n} & {\mbox{ if $s\in(0,\,1/2)$,}}\\
c\,|A|^{(n-1)/n} \log (|A|/|B|) & {\mbox{ if $s=1/2$,}}\\
c\,|A|^{(n-2s)/n} \,\big(|B|/|A|\big)^{1-2s} &{\mbox{ if $s \in(1/2,\,
1)$,}}
\end{matrix}
\right.\end{equation}
\item
if $|B| >c\,|A|$ then
\begin{equation}\label{gmt2}
L(A,D) \ge
c\,|A|^{(n-2s)/n} \,\big(|B|/|A|\big)^{-2s/n}
.\end{equation}
\end{itemize}\end{thm}

Theorems~\ref{EB}
and~\ref{DE} have, of course, physical relevance, since
they give optimal bounds on the energy of the limit
interface, and on the measure of the level sets
of the minimizers -- i.e., roughly speaking,
on the probability of finding a given phase
in a certain portion of the mediumm.

Also, due to the work of~\cite{CC}, density estimates
as the ones in Theorem~\ref{DE} are known to have
useful scaling properties and to play a crucial
role
in the geometric analysis of the level sets of the
rescaled minimizers,
especially in relation with the asymptotic interface. For instance,
we point out the following consequence of Theorem \ref{TH1} and Theorem~\ref{DE}.

\begin{corollary}\label{un}
Suppose that~$u$ is a global minimizer of $\EE$ in $\R^n$, i.e minimizes $\EE$ in any bounded domain $\Omega \subset \R^n$. Let
$$ u_\eps(x):=u\left(\frac{x}{\eps}\right).$$

Then
\begin{itemize}
\item[{(i)}] $u_\eps$ is a global minimizer for $\FF_\eps$,

\item[{(ii)}] $u_\eps$ converges, up to subsequences,
in~$L^1_{\rm
loc}(\R^n)$
to some~$u_\star =\chi_E-\chi_{\CC E}$ and $u_\star$ is a global minimizer of $\FF$ (see (\ref{FF1}), (\ref{FF2})),
\item[{(iii)}] given
any~$\theta\in(0,1)$, the set~$\{|u_\eps|\le\theta\}$
converges to $\partial E$ locally uniformly, that
is, for any~$R>0$ and any~$\delta>0$ there
exists~$\eps_o\in(0,1]$, possibly depending on~$R$ and~$\delta$,
such that, if~$\eps\in (0,\eps_o]$ then
\begin{equation}\label{IV}\{|u_\eps|\le\theta\}\cap B_R
\subseteq \bigcup_{p\in \partial E} B_\delta (p).\end{equation}
\end{itemize}
\end{corollary}

The minimizer $u$ above satisfies the 
Euler-Lagrange equation 
\begin{equation}\label{euler}
(-\Delta)^s u(x)+W'(u(x))=0,
\end{equation}
where $$(-\Delta)^s u(x):=\int_{\R^n}\frac{u(x)-u(y)}{|x-y|^{n+2s}}dy$$
and the integral is understood in the principal value sense. As usual,~$(-\Delta)^s$ is (up to a normalizing multiplicative
constant, depending on~$n$ and~$s$) the fractional power of
the positive operator~$-\Delta$.

Corollary \ref{un}
follows immediately: (i) from the scaling properties of 
$\EE$, (ii) from Theorem \ref{comp}, and (iii) is a consequence of the 
density estimates for the level sets of $u_\eps$ and the
$L^1_{\rm loc}$-convergence to $u_\star$ (see \cite{CC}
for further details).

The minimizing property of~$u_\star$ of~Corollary~\ref{un}(ii)
says that when~$s\in(0,\,1/2)$ the limit interface~$\partial E$ is
a nonlocal minimal surface in the setting of~\cite{CRS}, and when $s \in [
1/2, 1)$, $\partial E$ is a classical minimal surface.
This is interesting also because any regularity or rigidity
property proved for $\partial E$ may
reflect into similar ones for the minimizers of~$\FF_\eps$
(see, e.g.,~\cite{SA}).
In particular,~\eqref{euler}
may be seen as a semilinear equation driven by the fractional Laplacian.
Some rigidity properties for this kind of equations have been
recently obtained, for instance, in~\cite{CSM,SV}, but
many fundamental questions on this subject are still open.

The paper is organized as follows. In Section~\ref{GS} we prove 
Theorem~\ref{EB}. In Section~\ref{SS2} we prove Theorem~\ref{DE} 
by treating the cases~$s \in(0,1/2)$ and~$s\in[1/2,1)$ 
separately. 
Theorem~\ref{GMT}, together with a localized version of it 
(i.e., Proposition~\ref{GMTloc}), is proved in Section~\ref{GY}. 
Often in the proofs, when there is no possibility of confusion, we denote
the constants by~$C$ and $c$ although they may change from line to line.

\section{Proof of Theorem~\ref{EB}}\label{GS}

We use the following notation
$$u(A,B)=\int_A\int_B\frac{|u(x)-u(y)|^2}{|x-y|^{n+2s}}dxdy.$$

Since 
\begin{eqnarray*}
\EE(u,B_R)&=&\KK(u,B_R)+ \int_{B_R}W(u)\,dx
\\ &\le& \frac 12 u(B_{R+1},B_{R+1}) 
+ u(B_R, \CC B_{R+1}) + \int_{B_R}W(u)\,dx,\end{eqnarray*}
it suffices to bound each term
on the right by the quantity that appears in (\ref{hole}).

We define
\begin{equation}
\psi(x)=-1+2\min\{(|x|-R-1)^+,1\}
\end{equation}
so that $\psi=-1$ in $B_{R+1}$ and $\psi=1$ in $\CC B_{R+2}$.

First we show that $\EE(\psi,B_{R+2})$ satisfies the bounds in (\ref{hole}). Let $$ d(x):=\max\{R-|x|, \quad 1\}$$ and notice that

$$|\psi(x)-\psi(y)|\le \left\{
\begin{matrix}
2d(x)^{-1}|x-y| & {\mbox{ if $|x-y| < d(x)$,}}\\
2 &{\mbox{ if $|x-y|\ge d(x)$.}}
\end{matrix}
\right.$$

We obtain $$ \int_{\R^n}
\frac{
|\psi(x)-\psi(y)|^2
}{
|x-y|^{n+2s}}\,dy$$ $$\le \omega_{n-1}\int_0^{d(x)} \frac{(2r/d(x))^2}{r^{n+2s}}r^{n-1}dr + \omega_{n-1}\int_{d(x)}^\infty \frac{4}{r^{n+2s}}r^{n-1}dr \le C d(x)^{-2s}.$$

Now we integrate this inequality for all $x \in B_{R+2}$ and obtain that $\KK(\psi, B_{R+2})$ (therefore $\EE(\psi,B_{R+2})$) satisfies the energy bound of the theorem.

Let $v=\min\{u, \psi\}$ and denote $A:=\{v \le u\} \cap B_{R+2}.$ Clearly 
$B_{R+1} \subset A$ and $u=v$ in $\CC A$. We write that $u$ is a minimizer 
for $\EE$ in $B_{R+2}$, and therefore in~$A$:
\begin{equation}\label{497}\begin{split}&
\frac 12 u(A,A)+u(A, \CC A) + \int_A W(u)\,dx 
\\ &\qquad\le \frac 12 v(A,A)+v(A, 
\CC A) + \int_A W(v)\,dx.\end{split}\end{equation}
If $x \in A$ and $y\in \CC A$ then $v(x)=\psi(x) \le u(x)$, $v(y)=u(y) \le \psi(y)$ thus $$|v(x)-v(y)|\le \max\{|u(x)-u(y)|,|\psi(x)-\psi(y)|\},$$
which gives $$v(A, \CC A) \le u(A, \CC A)+ \psi(A, \CC A).$$
We use this in the energy inequality \eqref{497},
we simplify $u(A,\CC A)$ on both sides, and we obtain $$\frac 12 
u(A,A) + \int_A 
W(u)\,dx$$
 $$\le \frac 12 \psi(A,A)+\psi(A, \CC A) + \int_A W(v)\,dx = \EE(\psi, A) 
\le \EE(\psi, B_{R+2}).$$
Since $B_{R+1} \subset A$ we obtain the desired bounds for 
$u(B_{R+1},B_{R+1})$ and $\int_{B_R}W(u)\,dx$.

On the other hand $u(B_R, \CC B_{R+1})$ also satisfies a similar bound 
since $$\int_{\CC B_{R+1}}\frac{
|u(x)-u(y)|^2
}{
|x-y|^{n+2s}}\,dy \le C \int_{d(x)}^\infty r^{-1-2s} dr \le  C d(x)^{-2s} \quad \mbox{for all $x \in B_R$},$$
and then we integrate in $x\in B_R$.

\section{Proof of Theorem \ref{DE}}\label{SS2}

\subsection{Preliminary computations}

A minimizer $u$ of $\EE$ in $B_R$ with $R \ge 2$ satisfies the Euler-Lagrange 
equation in~\eqref{euler},
hence, $$\|u\|_{C^\alpha(B_1)} \le 
C(\|u\|_{L^\infty(\R^n)}+\| W'\|_{\R^n}) \le C,$$ 
for some small $\alpha >0$.

This shows that, by relabeling $\theta_1$, we can replace \eqref{MUCO1} by
\begin{equation}\label{MUCO}
\big| \{ u>\theta_1\}\cap B_{R_o}\big|\ge
\mu\end{equation}
for some constants~$R_o>0$ and~$\mu>0$,
possibly depending on $n$, $s$ and $W$ -- and, in fact,
the proof of Theorem~\ref{DE}
will make use only of~\eqref{MUCO}
rather than~\eqref{MUCO1}.

The strategy of the proof is, roughly speaking, to 
use the 
minimality of~$u$ in order to obtain an estimate of $|\{u>\theta_2\} \cap 
B_{2\rho}|$ in terms of $|\{u>\theta_2\} \cap B_\rho|$.
Then the conclusion will follow by iterating~\eqref{MUCO}. 

First, we construct the following useful barrier:

\begin{lemma}\label{tau}
Let~$n\ge1$. Given any~$\tau>0$, there exists~$C\ge1$, possibly
depending on~$n$, $s$ and $\tau$,
such that the following holds: for any $R\ge
C$, there exists a rotationally symmetric function
\begin{equation}\label{sopra}
w\in C\big( \R^n ,[-1+CR^{-2s},\,1]\big),
\end{equation}
with
\begin{equation}\label{sopra2}{\mbox{
$w=1$ in~$\CC B_R$,}}\end{equation} such that
\begin{equation}\label{al1} -(-\Delta)^s u(x)=
\int_{\R^n}\frac{w(y)-w(x)}{{|x-y|^{n+2s}}}\,dy
\le \tau\big(1+w(x)\big)\end{equation}
and
\begin{equation}\label{al2}
\frac1C \big( R+1-|x|\big)^{-2s}
\le 1+w(x)\le C \big( R+1-|x|\big)^{-2s}
\end{equation}
for any~$x\in B_R$.
\end{lemma}

\begin{proof} We fix a large $r\ge 1$, to be conveniently chosen
with respect to $R$ and $\tau$ in the sequel.

For $t\in(0,+\infty)$ and $x\in\R^n$, we define
\begin{eqnarray*}
g(t) &:=& t^{-2s},\\
h(t) &:=&\left\{
\begin{matrix}
\min\Big\{ 1,\, g(t)-g(r/2)-g'(r/2)\big( t-(r/2)\big)\Big\} \qquad
{\mbox{
if $t\le
r/2$,}}\\
0\qquad
{\mbox{ if $t\ge
r/2$,}}
\end{matrix}
\right. \\
v(x) &:=& \left\{ \begin{matrix} h(r-|x|) & {\mbox{ if $x\in B_r$,}}\\
1 & {\mbox{ if $x\in \CC B_r$.}}
\end{matrix}\right.
\end{eqnarray*} Such a~$v$, up to a proper rescaling, will
provide the existence of the desired function~$w$. To check this,
we first notice that
\begin{equation}
\label{3.5bis}
\begin{split}
&{\mbox{if $t\le r/2$ and $h(t)<1$ then}}\\
&\qquad h(t)=g(t)-g(r/2)-g'(r/2)\big( t-(r/2)\big)\\
&\qquad\qquad \ge g
(t)-g(r/2)-|g'(r/2)|\,(r/2)
\\&\qquad\qquad \ge t^{-2s}-16 r^{-2s}.\end{split}
\end{equation}
Moreover,
$v$ is continuous, radially symmetric, radially
nondecreasing and $0\le v\le 1$, due to the convexity
of $g$.

Also, we claim
that
\begin{equation}\label{giu}
{\mbox{
for any $x\in B_r$, $\| D^2 v\|_{ L^\infty( B_{ (r-|x|)/2 } (x)) }\le
2^8 (r-|x|)^{-2(1+s)}$.}}
\end{equation}
To prove \eqref{giu}, 
we observe that~$v=0$ in~$B_{r/2}$ and so~$D^2v=0$ in~$B_{r/2}$.
Then,
we take $y\in B_{ (r-|x|)/2 } (x) \cap (\CC B_{r/2})$
and we observe that
$$ |y|\le |y-x|+|x|\le \frac{r-|x|}2+|x|=
\frac{r+|x|}2,$$
hence
$$ r-|y|\ge r-\frac{r+|x|}2=\frac{r-|x|}2.$$
In particular,
$$ |r-|y||=r-|y|\le|y|$$
and, as a consequence,
\begin{eqnarray*} && |D^2 v(y)|\le\max\Big\{
|g''(r-|y|)| ,\,\frac{|g'(r-|y|)|}{
|r-|y||
}\Big\}\\
&&\qquad\le 2s(1+2s)\,{(r-|y|)^{-2(1+s)}}
\le \frac{2^8}{(r-|x|)^{-2(1+s)}},\end{eqnarray*}
proving \eqref{giu}.

{F}rom Lemma~6.15 
in~\cite{PSV} and~\eqref{giu}, we obtain that
\begin{equation}\label{giu1}
\begin{split}
&{\mbox{for any $x\in B_{r}$,}}\\
&\left| \int_{\R^n} \frac{
v(y)-v(x)
}{|x-y|^{n+2s}}\,dy\right|\\ &\qquad
\le C_0
\left[ \| D^2 v\|_{L^\infty(B_{(r-|x|)/2} (x))}
\Big(
(r-|x|)/2\Big)^{2(1-s)}
+ 2
\Big(
(r-|x|)/2\Big)^{-2s}
\right]
\\ &\qquad\le C_1 \Big[ (r-|x|)^{-4s}+(r-|x|)^{-2s}
\Big]
\\ &\qquad\le C_2 (r-|x|)^{-2s},
\end{split}
\end{equation}
for suitable $C_0$, $C_1$, $C_2>0$.

Now, we claim that
\begin{equation}\label{min 1}
\{ v<1\} \subseteq
B_{r-(1/2)}.
\end{equation}
Indeed, we take $x\in \{ v<1\}$ and
we define $t_x:=r-|x|$, so $h(t_x)<1$.
Hence either $t_x> r/2$
or $0<t_x\le r/2$ with $h(t_x)<1$.
In the first case, we would have that
$|x|\le r/2<r-1$ if $r$ is large enough and \eqref{min 1}
would hold, therefore we focus on the second case.
But then, recalling \eqref{3.5bis}, for large $r$,
\begin{eqnarray*} && 1>
t_x^{-2s} -{16}{r^{-2s}}
\ge t_x^{-2s}-(2^{2s}-1).
\end{eqnarray*}
That is, $t_x \ge 1/2$, proving \eqref{min 1}.

A straightforward consequence of \eqref{min 1} is that
\begin{equation}
\label{78H76} \| D^2 v\|_{L^\infty(\{ v<1\})} \le C_3,\end{equation}
for a suitable $C_3>0$.

Now, we set
$$ \Xi(x):=\left\{
\begin{matrix}
\nabla v(x) & {\mbox{ if $x\in \{ v<1\}$,}}\\
0 & {\mbox{ if $x\in \{ v=1\}$}}\\
\end{matrix}
\right.$$
and we
claim that
\begin{equation}
\label{78H77}
v(y)-v(x)-\Xi(x)\cdot(y-x)\le C_3 |x-y|^2,\end{equation}
for any~$x\in B_r$ and any~$y\in\R^n$.

To prove the claim above, we observe that, if~$x \in\{v=1\}$, then
the left hand side of~\eqref{78H77} is nonpositive, so~\eqref{78H77}
holds true. Also,~\eqref{78H77} follows from~\eqref{78H76} if
both~$x$ and~$y$ lie in~$\{v<1\}$, 
so it only remains
to prove~\eqref{78H77} when~$v(x)<1$ and~$v(y)=1$.
In such a case, we define~$v^\sharp$ to
be a smooth, radial extension of~$v$ outside~$\{v<1\}$
such that~$1\le v^\sharp \le2$ outside~$\{v<1\}$ and~$\| D^2
v^\sharp\|_{L^\infty(\{ v<1\})}
\le C_3$. Then,
\begin{eqnarray*}
&& v(y)-v(x)-\Xi(x)\cdot(y-x)= 1-v^\sharp(x)-\nabla v^\sharp(x)
\cdot(y-x)\\
&&\qquad\le v^\sharp (y)-v^\sharp(x)-\nabla v^\sharp(x)
\cdot(y-x)
\le C_3 |x-y|^2,
\end{eqnarray*}
proving~\eqref{78H77} in this case too.

Thus, thanks to~\eqref{78H77}, we
may use estimate~(6.47) in Lemma~6.14
of~\cite{PSV}
and obtain that
\begin{equation}\label{giu2}
\begin{split}
{\mbox{for any $x\in B_{r}$,}}\qquad
&\int_{\R^n} \frac{
v(y)-v(x)
}{|x-y|^{n+2s}}\,dy\,\le\, C_4
\end{split}\end{equation}
for a suitable $C_4>0$.

Now, we point out that
\begin{equation}\label{8.s}
\min\big\{ 1,t^{-2s}\big\}\le h(t)+16 r^{-2s}.
\end{equation}
Indeed, if $t\le r/2$, then \eqref{8.s}
is a consequence of \eqref{3.5bis}, while if $t>r/2$
we have that $t^{-2s}\le 8r^{-2s}<1$, which implies
\eqref{8.s}.

As a consequence of \eqref{8.s}, we have that
\begin{equation}\label{9.s}
\min\big\{ 1,(r-|x|)^{-2s}\big\}\le v(x)+16 r^{-2s}.
\end{equation}
{F}rom \eqref{giu1}, \eqref{giu2} and \eqref{9.s}, we conclude that
\begin{equation}\label{giu-e}
\begin{split}
{\mbox{for any $x\in B_{r}$,}}\qquad
&\int_{\R^n} \frac{
v(y)-v(x)
}{|x-y|^{n+2s}}\,dy \\&\qquad
\le C_5
\min\big\{ 1,(r-|x|)^{-2s}\big\} \\ &\qquad \le C_5
\big(
v(x)+16 r^{-2s}\big),\end{split}\end{equation}
for a suitable $C_5>0$.

Moreover, for any $t\in [0,\, r/2]$,
$$ h(t)\le g(t)+|g'(r/2)|\, (r/2)\le t^{-2s}+8r^{-2s}
$$
and so
\begin{equation}\label{giu-ee}
{\mbox{for any $x\in B_r$, }}
v(x)\le (r-|x|)^{-2s}+8r^{-2s}.
\end{equation}

Now we define
\begin{eqnarray*}
C_o & := & \Big( \frac{C_5}{\tau} \Big)^{1/2s},\\
\beta &:=& 32 r^{-2s}\\
w(x) &:=& (2-\beta) v\left(\frac{x}{C_o}\right)
+\beta-1.
\end{eqnarray*}
and we take $r:= R/C_o$.
Notice that $r$ is large if so is $R$, possibly
in dependence of $\tau$
(thus, from now on, the constants are
also allowed to depend on $\tau$).
Also, $w$ is radially non decreasing, and $w=1$ in $\CC B_R$.
In particular,
\begin{equation}\label{66}
1+w(x)\le 2.
\end{equation}
Moreover,
\begin{equation}\label{67}
{\mbox{
$w(x)=\beta-1$ for any $x\in B_{R/2}$}}
\end{equation}
and, from
\eqref{giu-ee}, for any $x\in B_R$,
$$ 1+w(x) \le 2v(x/C_o)+\beta \le 2C_o^{2s}
(R-|x|)^{-2s}+8r^{-2s} +\beta
.$$
That is, for a suitable $C_6>0$,
\begin{equation}\label{68}
1+w(x) \le C_6 (R-|x|)^{-2s}
\qquad{\mbox{
for any $x\in B_R\setminus
B_{R/2}$.}}\end{equation}
By \eqref{67} and \eqref{68}, we obtain that
\begin{equation}\label{69}
1+w(x) \le C_7 (R-|x|)^{-2s}
\qquad{\mbox{
for any $x\in B_R$,}}\end{equation}
for a suitable $C_7\ge 1$.

Now, we claim that
\begin{equation}\label{71}
1+w(x) \le \big(2^{2s}+2^{1-2s}\big) \,C_7\, (R+1-|x|)^{-2s}
\qquad{\mbox{
for any $x\in B_R$.}}\end{equation}
Indeed, if $|x|\le R-1$, we have that
$R-|x|\ge (R+1-|x|)/2$, therefore \eqref{71}
is a consequence of \eqref{69}. If, on the other
hand, $|x|\ge R-1$, we have that $R+1-|x|\le 2$,
thus we
use \eqref{66} to obtain
$$ 1+w(x)\le 2=2^{1-2s} 2^{-2s}
\le 2^{1-2s}
(R+1-|x|)^{-2s},$$
which gives \eqref{71}.

Then, \eqref{71} implies the upper bound in~\eqref{al2}, and the lower
bound may be obtained analogously, using~\eqref{3.5bis}.

Moreover, recalling \eqref{giu-e},
for any $x\in B_{R}$,
\begin{eqnarray*}
&&\int_{\R^n} \frac{
w(y)-w(x)
}{|x-y|^{n+2s}}\,dy
= (2-\beta)\, C_o^{-2s}
\int_{\R^n} \frac{
v(y)-v(x/C_o)
}{|(x/C_o)-y|^{n+2s}}\,dy
\\ &&\qquad
\le (2-\beta)\, C_o^{-2s} C_5 \Big( v(x/C_o)+16 r^{-2s}\Big)\\ &&\qquad
\le C_o^{-2s} C_5 \Big( (2-\beta) v(x/C_o)+32 r^{-2s}\Big)\\ &&\qquad
=\tau \big(1 +w(x)\big).
\end{eqnarray*}
This proves \eqref{al1}
and it completes the proof of Lemma~\ref{tau}.
\end{proof}

Now, we give an elementary, general estimate:

\begin{lemma}\label{lemma indu}
Let~$\sigma$, $\mu\in(0,+\infty)$, $\nu\in(\sigma,+\infty)$
and~$\gamma$, $R_o$, 
$C\in(1,+\infty)$.

Let~$V:(0,+\infty)\rightarrow(0,+\infty)$ be a
nondecreasing function. For any~$r\in [R_o,+\infty)$, let
$$ \alpha(r):=\min\left\{1,\,\frac{\log V(r)}{\log r}
\right\}.$$
Suppose that
\begin{equation}\label{ind 1}
V(R_o)\ge\mu
\end{equation}
and
\begin{equation}\label{ind 2}
r^\sigma\,\alpha(r)\,V(r)^{(\nu-\sigma)/\nu}\le C V(\gamma r),\qquad
{\mbox{ for any $r\in[R_o,+\infty)$.}}
\end{equation}
Then, there exist~$c\in(0,1)$ and~$R_\star\in[R_o,+\infty)$, possibly 
depending
on~$\mu$, $\nu$, $\gamma$, $R_o$ and~$C$, such that
$$ V(r)\ge cr^\nu,\qquad
{\mbox{ for any $r\in[R_\star,+\infty)$.}}$$
\end{lemma}

\begin{proof} Let~$j_1$ be the smallest
natural number for which~$\gamma^{j_1}\ge R_o$. Notice that
such a definition
is well posed since~$\gamma>1$.

Let
\begin{equation}
\label{ind 3} c:=\min\left\{\frac{\mu}{\gamma^{\nu 
j_1}},\,\left(
\frac1{C\gamma^\nu}\right)^{\nu/\sigma},\,
\left(\frac\nu{2C\gamma^\nu}\right)^{\nu/\sigma}
\right\}.\end{equation}

Let~$j_2$ be the smallest integer for which
\begin{equation}\label{ind 4}\frac{|\log c|}{j_2 \log\gamma}\le
\frac\nu2.
\end{equation}

Let~$j_\star:= j_1+j_2$. For any~$j\in\N\cap [j_\star,+\infty)$, we
define~$v_j:= V(\gamma^j)$. We claim that
\begin{equation}\label{ind 5}
v_j\ge c \gamma^{\nu j},\qquad
{\mbox{ for any $r\in[R_o,+\infty)$.}}
\end{equation}
The proof of~\eqref{ind 5} is by induction over~$j$.
First of all,
$$ v_{j_\star}\ge V(\gamma^{j_1})\ge V(R_o)\ge\mu\ge c\gamma^{\nu 
j_1},$$
thanks to~\eqref{ind 3}.

Then, we suppose that~\eqref{ind 5} holds for some~$j\ge j_\star$
and we prove it for~$j+1$, via the following argument.

Since we assumed that~\eqref{ind 5} holds~$j$,
\begin{eqnarray*}
\alpha(\gamma^j) &=& \min\left\{1,\,\frac{\log (v_j)}{\log \gamma^j}
\right\}\\
&\ge&
\min\left\{1,\,\frac{\log (c\gamma^{\nu j})}{\log \gamma^j}
\right\}
\\ &=&
\min\left\{1,\,\nu-\frac{|\log c|}{j \log \gamma}
\right\}
\\ &\ge& \min\left\{1,\,\frac{\nu}2
\right\},
\end{eqnarray*}
thanks to~\eqref{ind 4}.

Therefore, using~\eqref{ind 2} with~$r:=\gamma^j$
and the assumption that~\eqref{ind 5} holds for~$j$,
we conclude that
\begin{eqnarray*}
v_{j+1} &=& V(\gamma^{j+1})
\\ &\ge&
\frac{\gamma^{\sigma j}}{C}\,\alpha(\gamma^j)\,v_j^{(\nu-\sigma)/\nu}
\\ &\ge&
\min\left\{\frac1C,\,\frac{\nu}{2C}
\right\} c^{(\nu-\sigma)/\nu} \gamma^{\nu j}
.\end{eqnarray*}
Recalling~\eqref{ind 3}, we see that this last quantity
is greater or equal than~$c \gamma^{(j+1)\nu}$, thus
completing the induction argument which proves~\eqref{ind 5}.

{F}rom~\eqref{ind 5}, the desired result plainly follows.
\end{proof}

\begin{rem}\label{R indu}
{\rm
In the sequel, we will use Lemma~\ref{lemma indu}
with~$V(R):= | \{ u>\theta_\star\}\cap
B_R|$, with~$\theta_\star\le\theta_1$.
In this way, condition~\eqref{ind 1}
is warranted by~\eqref{MUCO}.
}\end{rem}

Now, we make some useful computations, valid for any~$s\in(0,1)$.
We fix~$K\ge 2(R_o+1)$, to be taken suitably
large in the sequel, where~$R_o$ is the one given by
the statement of Theorem~\ref{DE}, and~$R>2K$.
Given any measurable~$w:\R^n\rightarrow[-1,1]$ such that
\begin{equation}\label{la34}
{\mbox{$w= 1$ in $\CC B_R$,}}
\end{equation}
we define
\begin{equation}\label{3.27bis}
v(x):=\min\{u(x),w(x)\}\end{equation}
and~$D:=(\R^n\times\R^n)\setminus
(\CC B_R\times \CC B_R)$.
Notice that
$$\KK(u;B_R)=
\frac{1}{2}\iint_D
\frac{
|u(x)-u(y)|^2
}{
|x-y|^{n+2s}}\,dx\,dy$$
and that
\begin{equation}\label{agree}
{\mbox{
$v=u$ in~$\CC B_R$.}}\end{equation}
So, by a simple algebraic computation, we have that
\begin{eqnarray*}&&
\KK(u-v;B_R)+\KK(v;B_R)-\KK(u;B_R)\\
&&\qquad=
-\iint_D \frac{\big( (u-v)(x)-(u-v)(y)\big)\,\big(
v(y)-v(x)\big)}{|x-y|^{n+2s}}\,dx\,dy\\
&&\qquad =
-\iint_{\R^n\times\R^n} \frac{\big( (u-v)(x)-(u-v)(y)\big)\,\big(
v(y)-v(x)\big)}{|x-y|^{n+2s}}\,dx\,dy
\\
&&\qquad = -2
\iint_{\R^n\times\R^n} \frac{\big( u(x)-v(x)\big)\,\big(
v(y)-v(x)\big)}{|x-y|^{n+2s}}\,dx\,dy
\\
&&\qquad =
2
\int_{ \{ u>v=w\}} \big( u(x)-v(x)\big)\, \left[
\int_{\R^n}\frac{
v(y)-w(x)}{|x-y|^{n+2s}}\,dy\right]\,dx
\\
&&\qquad \le
2
\int_{ B_R\cap \{ u>w\}} \big( u(x)-v(x)\big)\, \left[
\int_{\R^n}\frac{
w(y)-w(x)}{|x-y|^{n+2s}}\,dy\right]\,dx
\\
&&\qquad =-
2
\int_{ B_R\cap \{ u>w\}} \big( u(x)-v(x)\big)\, (-\Delta)^s w(x)
\,dx
.\end{eqnarray*}
As a consequence, using once more~\eqref{agree} and
the minimality of~$u$, we conclude that
\begin{equation}\label{la3}
\begin{split}
&\KK(u-v;B_R)\\
&\qquad\le \EE_R(u)-\EE_R(v)
+\int_{B_R} W(v)-W(u)\,dx\\
&\qquad\quad-
2
\int_{ B_R\cap\{ u>w\}} \big( u(x)-v(x)\big)\, (-\Delta)^s w(x)
\,dx
\\ &\qquad\le
\int_{B_R\cap \{ u>w\}} W(w)-W(u)\,dx
\\
&\qquad\quad-
2
\int_{ B_R\cap\{ u>w\}} \big( u(x)-w(x)\big)\, (-\Delta)^s w(x)
\,dx
.
\end{split}
\end{equation}
Now, we fix~$\theta_1$, $\theta_2 \in(-1,1)$
as in the statement of
Theorem~\ref{DE}, and we take
\begin{equation}\label{TS}
\theta_\star\le\min\{\theta_1,\,\theta_2,\, -1+c\},
\end{equation}
with~$c$ as in~\eqref{grow}. We define
\begin{equation}\label{CHS}
A(R)\,:=\,c\,\int_{ B_R \cap \{w<u\le\theta_\star\}}
(u-w)^2\,dx.
\end{equation}
{F}rom the behavior of~$W$ near its minima (see~\eqref{grow}), we
deduce that
\begin{equation*}\begin{split}
&\int_{B_R\cap\{u>w\}} W(w)-W(u)\,dx\\
&\quad\le
-c\int_{B_R\cap \{w<u\le \theta_\star\}} (1+w)(u-w)\,dx
\\ &\qquad+\frac{1}{c}
\int_{B_R\cap \{u>\max\{w,\theta_\star\}\} } (1+w)\,dx
\,-\,A(R)
.
\end{split}\end{equation*}
This and~\eqref{la3}
give that
\begin{equation}\label{la4}
\begin{split}
\KK(u-v;B_R)&\le
-c\int_{B_R\cap \{w<u\le\theta_\star\}} (1+w)(u-w)\,dx
\\
&\qquad
+\frac{1}{c}
\int_{B_R\cap \{u>\max\{w,\theta_\star\}\} } (1+w)\,dx
\\
&\qquad-
2
\int_{ B_R\cap \{ u>w\}} \big( u(x)-w(x)\big)\, (-\Delta)^s w(x)
\,dx
\\ &\qquad-\,A(R).
\end{split}
\end{equation}
While~\eqref{la4} is valid for any~$w$
satisfying~\eqref{la34}, we now choose~$w$ in a convenient way.
That is, we define~$\tau:= c/4$ and we take~$w$ to be
the function constructed in Lemma~\ref{tau}.
With this choice,~\eqref{la4} and
Lemma~\ref{tau} give that
\begin{equation}\label{la5}
\begin{split}
&\KK(u-v;B_R)
+\frac{c}2\int_{B_R\cap \{w<u\le \theta_\star\}} (1+w)(u-w)\,dx
\\ &\le
-\frac{c}2\int_{B_R\cap \{w<u\le \theta_\star\}} (1+w)(u-w)\,dx
\\
&\qquad
+\frac{1}{c}
\int_{B_R\cap \{u>\max\{w,\theta_\star\}\} } (1+w)\,dx
\\&\qquad+
2\tau
\int_{ B_R\cap \{ u>w\}} \big( u(x)-w(x)\big)\,
\big( 1+w(x)\big)\,dx\,-\,A(R)
\\&\le
\frac{1}{c}
\int_{B_R\cap\{u> \max\{w,\theta_\star\} \}} (1+w)\,dx
\,-\,A(R)
\\ &\le
C \int_{B_R\cap \{ u>\theta_\star\} } (R+1-|x|)^{-2s}\,dx
\,-\,A(R)
,
\end{split}
\end{equation}
for a suitable~$C>0$.

Now, we define
\begin{equation}\label{NOTV}
V(R):=\big| \{ u>\theta_\star\}\cap
B_R\big|.\end{equation}
Hence, using the Coarea formula, we deduce from~\eqref{la5}
that
\begin{equation}\label{la6}
\begin{split}&
\!\! A(R)+\KK(u-v;B_R)
+\frac{c}2\int_{B_R\cap \{w<u\le \theta_\star\}} (1+w)(u-w)\,dx
\\
&\le
C \int_{0}^R (R+1-t)^{-2s}
\left( \int_{\partial B_t}
\chi_{\{ u>\theta_\star\}}(x)\,d{\mathcal{H}}^{n-1}(x)\right)\,dt
\\ &=C \int_{0}^R (R+1-t)^{-2s} V'(t)\,dt.
\end{split}
\end{equation}

\subsection{Completion of the proof of~Theorem \ref{DE}}

Now we use the estimates of Theorem~\ref{GMT}. For
convenience, the proof
of Theorem~\ref{GMT} itself is postponed to Section~\ref{GY}.


Given a
measurable set~$A\subseteq\R^n$,
we define
\begin{equation}\label{3.36bis} \ell (A):=
\left\{
\begin{matrix}
|A|^{(1-2s)/n} & {\mbox{ if $s\in(0,\,1/2)$,}}\\    
\log|A| &{\mbox{ if $s=1/2$,}}\\
1 & {\mbox{ if $s\in (1/2,\,1)$.}}
\end{matrix}
\right. \end{equation}
Given~$s\in(1/2,1)$ and~$\alpha\ge0$, we notice that the map
$$ (0,+\infty)\ni t\mapsto\alpha t^{1-2s}+t$$
has minimum at~$t=\big(\alpha(2s-1)\big)^{1/(2s)}$ and therefore
\begin{equation}\label{77676} \inf_{t\in (0,+\infty)}
\alpha t^{1-2s}+t\ge c_1 \alpha^{1/(2s)},
\end{equation}
for a suitable~$c_1>0$, as long as~$s\in(1/2,1)$.

Also, the map
$$ (0,+\infty)\ni t\mapsto\alpha^{(n-1)/n}\log\frac{\alpha}{t}+t$$
has minimum at~$t=\alpha^{(n-1)/n}$, so
\begin{equation}\label{77676.1}
\inf_{t\in(0,+\infty)}
\alpha^{(n-1)/n}\log\frac{\alpha}{t}+t\ge
\frac{\alpha^{(n-1)/n}}n \log\alpha.
\end{equation}
Moreover, if, given~$\kappa>0$, we consider the map
$$ [\kappa,+\infty)\ni t\mapsto
\Phi(t):=\frac{t^{1/n}}{1+|\log t|},$$
we have that~$\Phi(t)>0$ for any~$t\in [\kappa,+\infty)$
and
$$\lim_{t\rightarrow+\infty}\Phi(t)=+\infty,$$
therefore
\begin{equation}\label{77676.2}
i_\kappa := \inf_{t\in [\kappa,+\infty)} \Phi(t)>0.\end{equation}
Now, we claim that, if~$A$ and~$B$
are disjoint measurable subsets of~$\R^n$ and~$D:=\CC(A\cup
B)$, with~$|A|\ge \kappa>0$,
then there exists~$c_0>0$, possibly depending on~$\kappa$
such that
\begin{equation}\label{sum}
L(A,D)+|B|\,\ge \,c_0\,|A|^{(n-1)/n} \,\ell (A).
\end{equation}
To prove~\eqref{sum}, we take~$c$ as\footnote{Of course,
no confusion should arise with the~$c$ in~\eqref{grow}
which was previously used in~\eqref{TS} -- in any case,
one could just take~$c$ to be the smallest
of these two constants to make them equal.}
in Theorem~\ref{GMT}
and we distinguish two cases. First, if~$|B|>c|A|$, we
use~\eqref{77676.2} to see that
\begin{eqnarray*}&&
L(A,D)+|B|\ge |B|> c|A|
\\ && =\left\{\begin{matrix}
c|A|^{2s/n} |A|^{(n-2s)/n}\ge c\kappa^{2s/n}\,
|A|^{(n-2s)/n},\\
\qquad{\mbox{ if $s\in(0,1/2)$,}}\\
\, \\
c|A|^{1/n} |A|^{(n-1)/n}\ge c\,i_\kappa\,
(1+|\log|A||)   
|A|^{(n-1)/n},\\
\qquad{\mbox{ if $s\in[1/2,1)$,}}\end{matrix}
\right.\end{eqnarray*}
which gives~\eqref{sum}.
Therefore, we may suppose that~$|B|\le c|A|$ and
use~\eqref{gmt1},
\eqref{77676} and~\eqref{77676.1}
to conclude that
\begin{eqnarray*}
\frac{1}{c} \Big( L(A,D)+|B|\Big)& \ge& \left\{
\begin{matrix}
|A|^{(n-2s)/n} +|B| & {\mbox{ if $s\in(0,\,1/2)$,}}\\
|A|^{(n-1)/n} \log (|A|/|B|)+|B| & {\mbox{ if $s=1/2$,}}\\
|A|^{2s(n-1)/n} |B|^{1-2s} +|B|&{\mbox{ if $s \in(1/2,\,
1)$,}}
\end{matrix}
\right.
\\& \ge& \left\{
\begin{matrix}
|A|^{(n-2s)/n} & {\mbox{ if $s\in(0,\,1/2)$,}}\\    
\Big(
|A|^{(n-1)/n} \log |A|\Big)/n & {\mbox{ if $s=1/2$,}}\\
c_1 \Big(|A|^{2s(n-1)/n}\Big)^{1/(2s)} &{\mbox{ if $s \in(1/2,\,
1)$,}}
\end{matrix}
\right.\\ &\ge& c_2\,|A|^{(n-1)/n}\,\ell(A),
\end{eqnarray*}
for some~$c_2>0$, proving~\eqref{sum}.

Now, we take a free parameter~$K>1$, that will
be chosen conveniently large in what follows. The radius~$R$
of Theorem~\ref{DE}
will be taken larger than~$K$.
We observe that,
by~\eqref{al2},   
\begin{equation}\label{RON}{\mbox{
$w\le
-1+C (K+1)^{-2s}
<-1+\displaystyle\frac{1+\theta_\star}2$ in~$B_{R-K}$,}}\end{equation}
as long as~$K$
is large enough possibly in dependence of~$\theta_\star$
which was fixed in~\eqref{TS},
and so
\begin{equation} \label{CU} a_R:=B_{R}\cap \left\{u-w\ge
\frac{1+\theta_\star}{4}\right\}
\supseteq B_{R-K}\cap \{ u>\theta_\star\}.\end{equation}
By~\eqref{MUCO}, \eqref{TS} and~\eqref{CU}, when $R$ is large
\begin{equation*}
|a_R|\ge|\{u>\theta_1\}\cap B_{R_o}|\ge \mu.
\end{equation*}
As a consequence, we may apply~\eqref{sum}
with
\begin{eqnarray*} && \kappa:=\mu,\\&& A:=a_R,\\
&& B:=b_R:=
B_{R}\cap
\left\{ \frac{1+\theta_\star}{8}<
u-w
<\frac{1+\theta_\star}{4}
\right\}\\ {\mbox{and }}&&
D:=d_R:=\CC (A\cup B)=
\CC B_R\cup\left(
B_R\cap \left\{
u-w\le\frac{1+\theta_\star}{8}
\right\} \right)
.\end{eqnarray*}
We obtain that
\begin{equation*}\begin{split}&
L\left(
a_R,d_R
\right)+\left|b_R\right|\ge
\,c_3\, \left|
a_R \right|^{(n-1)/n}\,
\ell
\left(
a_R\right)
\end{split}
\end{equation*}
for a suitable~$c_3>0$, possibly depending on~$\mu$.

Accordingly, recalling
the notation in~\eqref{NOTV}, \eqref{3.36bis} and~\eqref{CU},
\begin{equation}\label{5653jcjama23}
\begin{split}&
L\left(
a_R,d_R
\right)+\left| b_R\right|\,\ge\, c_3 \,V(R-K)^{(n-1)/n}\,\ell_{R-K},
\end{split}\end{equation}
where
\begin{equation}\label{3.43bis} \ell_R :=\ell\Big(
B_R\cap \{u>\theta_\star\}\Big)=\left\{
\begin{matrix}
(V(R))^{(1-2s)/n}& {\mbox{ if $s\in (0,\,1/2)$,}}\\
\log (V(R)) & {\mbox{ if $s=1/2$,}}\\
1 & {\mbox{ if $s\in (1/2,\,1)$.}}
\end{matrix}
\right.\end{equation}
Notice that
\begin{equation}\label{nod}{\mbox{the map $R\mapsto \ell_R$
is nondecreasing.}}\end{equation}
Now, we recall~\eqref{3.27bis} and~\eqref{agree}
to see that
\begin{equation*}
\begin{split}
&{\mbox{if }} y\in d_R
{\mbox{ then }}
(u-v)(y) \le
\frac{1+\theta_\star}{8}
.\end{split}\end{equation*}
Therefore,
\begin{equation*}
\begin{split}
&{\mbox{if }} y\in d_R
{\mbox{ and }} x\in a_R
\,
{\mbox{ then }}
\\
&\qquad |(u-v)(x)-(u-v)(y)|\ge
(u-v)(x)-(u-v)(y)\\ &\qquad\qquad\ge
\frac{1+\theta_\star}{4}-\frac{1+\theta_\star}{8}
=\frac{1+\theta_\star}{8}
.\end{split}\end{equation*}
Recalling~\eqref{1.8bis}, this implies that
\begin{equation}\label{67dfu556}
\begin{split}& 2\,\KK(u-v;B_R)\\
&\qquad\ge
\int_{a_R}\int_{d_R}
\frac{|(u-v)(x)-(u-v)(y)|^2}{|x-y|^{n+2s}}\,dx\,dy
\\ &\qquad\ge \left( \frac{1+\theta_\star}{8}
\right)^2 L\left(
a_R,d_R\right).\end{split}\end{equation}
Furthermore, by~\eqref{RON},
\begin{equation}\label{onr}
b_{R-K}=B_{R-K}\cap
\left\{ \frac{1+\theta_\star}{8}<
u-w
<\frac{1+\theta_\star}{4}
\right\} \subseteq B_{R-K}\cap \{ u\le\theta_\star\}.\end{equation}
Also, by~\eqref{CHS},
\begin{equation}\label{es-ar}\begin{split}
& A(R)=c\int_{B_R\cap \{w<u\le\theta_\star\} } (u-w)^2\,dx
\\ &\qquad\ge c
\int_{ b_{R} \cap \{u\le\theta_\star\}} (u-w)^2\,dx
\\ &\qquad\ge
c_4 \,\left|b_{R} \cap \{u\le\theta_\star\}\right|,
\end{split}\end{equation}
for a suitable~$c_4\in(0,1)$.

Now, we observe that, if~$t\in[R-K,R]$,
$$ (R+1-t)^{-2s}\ge (1+K)^{-2s}$$
and therefore, recalling the notation in~\eqref{NOTV},
\begin{eqnarray*}
&&V(R)-V(R-K)=\int_{R-K}^R V'(t)\,dt\\
&&\qquad\le \frac{1}{c_5} \int_{R-K}^R (R+1-t)^{-2s}V'(t)\,dt
\le \frac{1}{c_5} \int_{0}^R (R+1-t)^{-2s}V'(t)\,dt,
\end{eqnarray*}
for a suitable~$c_5\in(0,1)$, depending on~$K$,
that is now fixed once and for all.
Therefore, exploiting~\eqref{onr},
\begin{eqnarray*}
|b_R|&\le& |b_R \cap \{u\le\theta_\star\}|+
|b_R \cap \{u > \theta_\star\}|
\\ &=& \big|b_R \cap \{u\le\theta_\star\}\big|+
\big|(b_R\setminus b_{R-K})\cap \{u > \theta_\star\}\big|
\\ &\le&
\big|b_R \cap \{u\le\theta_\star\}\big|+
\big|(B_R\setminus B_{R-K})\cap \{u > \theta_\star\}\big|
\\ &=&
|b_R \cap \{u\le\theta_\star\}|+
V(R)-V(R-K)\\ &\le&
|b_R \cap \{u\le\theta_\star\}|+
\frac{1}{c_5} \int_{0}^R (R+1-t)^{-2s}V'(t)\,dt.
\end{eqnarray*}
With this, we can write~\eqref{es-ar} as
\begin{eqnarray*} && A(R)+
\frac{c}2\int_{B_R\cap \{w<u\le \theta_\star\}} (1+w)(u-w)\,dx
\ge A(R)\\ &&\qquad\ge
c_4 \left|b_{R} \right|
-
\frac{c_4}{c_5} \int_{0}^R (R+1-t)^{-2s}V'(t)\,dt.
\end{eqnarray*}
The latter estimate, \eqref{la6} and~\eqref{67dfu556}
give that
\begin{eqnarray*}
&& C \int_{0}^R (R+1-t)^{-2s} V'(t)\,dt\\
&\ge& \KK(u-v;B_R)+A(R)+
\frac{c}2\int_{B_R\cap \{w<u\le \theta_\star\}} (1+w)(u-w)\,dx
\\ &\ge& c_6 \Big(
L\left(a_R,d_R
\right)
+\left|b_R\right|\,\Big)
-\frac{c_4}{c_5} \int_{0}^R (R+1-t)^{-2s}V'(t)\,dt
,\end{eqnarray*}
for suitable~$C\in(1,+\infty)$ and~$c_6\in(0,1)$, possibly depending
on~$\theta_\star$, that was fixed in~\eqref{TS}.

Therefore, taking the last term on the other
side, recalling~\eqref{5653jcjama23}
and possibly renaming~$C\ge1$, we conclude that\footnote{The
reader may observe in~\eqref{F5}
the structurally different iteration between the
cases~$s\in(0,\,1/2)$
and~$s\in[1/2,\,1)$, which is encoded in~$\ell_{R-K}$,
according to~\eqref{3.43bis}.}
\begin{equation}\label{F5}
\begin{split}&
C \int_{0}^R (R+1-t)^{-2s} V'(t)\,dt\,\ge\,
\ell_{R-K} V(R-K)^{(n-1)/n}
.\end{split}\end{equation}
Now, we notice that, if~$\rho\ge1$,
\begin{equation}\label{dhjdeu2}
\int_t^{(3/2)\rho}(R+1-t)^{-2s} \,dR \,\le\,\left\{ \begin{matrix}
(4^{1-2s} /(1-2s))\rho^{1-2s}  & {\mbox{ if $s\in(0,\,1/2)$,}}\\    
\log((5/2)\rho) &{\mbox{ if $s=1/2$,}} \\ 1/(2s-1) & {\mbox{ if
$s\in(1/2,\,1)$.}} \end{matrix} \right.\end{equation}

Therefore, we make use of~\eqref{nod} and~\eqref{dhjdeu2}
in order to integrate~\eqref{F5} in~$R\in
[\rho,\,(3/2)\rho]$, with~$\rho\ge 2K$, and we obtain that
\begin{equation*} \begin{split}
& \rho\,
\ell_{\rho-K}\,V(\rho-K)^{(n-1)/n} \\
&\qquad\le C \int_\rho^{(3/2)\rho}\left( \int_{0}^R (R+1-t)^{-2s}
V'(t)\,dt \right)\,dR \\ &\qquad \le C \int_0^{(3/2) \rho}\left( \int_{
t
}^{(3/2)\rho} (R+1-t)^{-2s} \,dR
\right)\,V'(t)\,dt
\\ &\qquad \le \,\left\{
\begin{matrix}
C' \rho^{1-2s} V((3/2)\rho) & {\mbox{ if $s\in(0,\,1/2)$,}}\\    
C' \,(\log\rho)\, V((3/2)\rho) & {\mbox{ if $s=1/2$,}}\\
C' V((3/2)\rho) &{\mbox{ if $s\in(1/2,\,1)$.}}
\end{matrix}
\right.
\end{split}\end{equation*}
for some $C$, $C'\ge1$.

That is, for large $r$, recalling~\eqref{3.43bis}
and possibly renaming $C>1$,
\begin{equation}\label{7d7dhjjj0}
\left\{\begin{matrix}
r^{2s} V(r)^{(n-2s)/n}\le CV(2r) & {\mbox{if $s\in(0,\,1/2)$,}}\\
r\displaystyle\frac{\log V(r)}{\log r}V(r)^{(n-1)/n}\le CV(2r) 
& {\mbox{if $s=1/2$,}}\\
rV(r)^{(n-1)/n}\le CV(2r) & {\mbox{if $s\in(1/2,\,1)$.}}
\end{matrix}\right.
\end{equation}
By~\eqref{7d7dhjjj0}
and Lemma~\ref{lemma indu}, applied here with~$\sigma:=2s$ when~$s\in
(0,\,1/2)$, or~$\sigma:=1$ when~$s\in[1/2,\,1)$,
and~$\gamma:=2$ (recall also Remark~\ref{R indu}),
we obtain 
that~$V(R)\ge c_o R^n$ for large~$R$,
for a suitable~$c_o\in(0,1)$.

Therefore, recalling also~\eqref{TS},
\begin{equation}\label{PP1}\begin{split}
& \big| \{ u>\theta_2 \}\cap B_R\big|+
\big| \{ \theta_\star<u\le \theta_2
\}\cap B_R\big|\\
&\qquad=\big| \{ u>\theta_\star\}\cap B_R\big| = V(R)\ge c_o
R^n\end{split}
\end{equation}
for large~$R$.
On the other hand, by~\eqref{hole},
\begin{equation}\label{PP2}\begin{split}
&\frac{\overline{C} R^{n}}{\log R}\ge \EE(u,B_R)\ge 
\int_{ \{ \theta_\star<u\le \theta_2
\}\cap B_R } W(u(x))\,dx
\\&\qquad \ge\,\inf_{r\in[\theta_\star,
\theta_2]} W(r)\;
\big| \{ \theta_\star<u\le \theta_2
\}\cap B_R\big|
.\end{split}\end{equation}
By~\eqref{PP1} and~\eqref{PP2},
we obtain that~\eqref{hard} holds true,
and this completes the proof of Theorem~\ref{DE}.

\section{Proof of Theorem~\ref{GMT}}\label{GY}

Given $\xi\in {\rm S}^{n-1}$, we denote by $\pi_\xi$ the hyperplane
normal to $\xi$ passing through the origin, namely
$$ \pi_\xi := \big\{ x\in \R^n {\mbox{ s.t. }} \xi\cdot  x=0\}.$$
Given $\Omega\subseteq\R^n$,
we consider the projection of $\Omega$ along $\xi$, i.e.
$$ \Pi_\xi(\Omega):=\big\{
p\in \pi_\xi {\mbox{ s.t. there exists $t\in\R$ for which }}
p +t\xi\in \Omega
\big\}.$$
Next result relates the $n$-dimensional
measure of $\Omega$ with the largest possible
$(n-1)$-dimensional measure of $\Pi_\xi(\Omega)$
(i.e., pictorically, the measure of an object in a room
with the measure of its shadows on the walls and on the floor).

\begin{lemma}\label{proj}
Let $\Omega$ be a measurable subset of $\R^n$.
Then,
\begin{itemize}
\item[{(i)}] $|\Omega|^{n-1}\le
|\Pi_{e_1}(\Omega)|\dots|\Pi_{e_n}(\Omega)|$,
\item[{(ii)}] there exists $k\in \{1,\dots,n\}$
for which
$$ |\Pi_{e_k} (\Omega)|\,\ge \,|\Omega|^{(n-1)/n}.$$
\end{itemize}\end{lemma}

\begin{proof} First of all, we use the generalized H\"older
inequality (see, e.g., page~623 of~\cite{Evans}) to observe that,
if~$\psi_o\ge0$ and~$\psi_1,\dots,\psi_{n-1}
\in L^1\big(\R,\,[0,+\infty]\big)$, then
\begin{equation}\label{oh}
\begin{split}
& \int_\R \Big( \psi_o \psi_1(t)\dots\psi_{n-1}(t)\Big)^{1/(n-1)}\,dt
\\&\qquad=\psi_o^{1/(n-1)}\,\int_\R \Big(
\psi_1(t)\dots\psi_{n-1}(t)\Big)^{1/(n-1)}\,dt
\\ &\qquad\le\psi_o^{1/(n-1)}
\left(\int_\R\psi_1(t)\,dt\right)^{1/(n-1)}
\dots\left(\int_\R\psi_{n-1}(t)\,dt\right)^{1/(n-1)}
\\ &\qquad=\left(\psi_o
\int_\R\psi_1(t)\,dt
\dots\int_\R\psi_{n-1}(t)\,dt\right)^{1/(n-1)}.
\end{split}\end{equation}
Now, we introduce some notation. Given~$x=(x_1,\dots,
x_n)\in\R^n$, for any~$i\in\{1,\dots,n\}$
we define
$$\hat x_i:=x-x_ie_i=(x_1,\dots,x_{i-1},0,x_{i+1},\dots,x_n).$$
Also, for any~$i\le k\in\{1,\dots,n\}$, we set
$$ \hat x_{i;k}:=
(x_1,\dots,x_{i-1},x_{i+1},\dots,x_k).$$
Notice that~$\hat x_i\in\R^n$, $\hat x_{i;k}\in\R^{k-1}$
and the above notation means that~$\hat x_{k;k}:=(x_1,\dots,x_{k-1})$.
We also stress the fact that
\begin{equation}\label{oh2}
{\mbox{$\hat x_i$ does not depend on~$x_i$.}}
\end{equation}
For short, we also
set~$\chi:=\chi_\Omega$ and~$\chi_i:=\chi_{\Pi_{e_i}(\Omega)}$.
Then, we have that
\begin{equation}\label{oh3}
\int_{\R^{n-1}} \chi_i(\hat x_i)
d\hat x_{i;n}
=|\Pi_{e_i}(\Omega)|.
\end{equation}
Now, we observe that
$$ \chi(x)\le \chi_1(\hat x_1)\dots\chi_n(\hat x_n)$$
and so
$$ \chi(x)=\Big(\chi(x)\Big)^{1/(n-1)}
\le \Big(\chi_1(\hat x_1)\dots\chi_n(\hat
x_n)\Big)^{1/(n-1)}.$$
Hence, integrating in~$dx_1$ and using~\eqref{oh}
and~\eqref{oh2},
$$ \int_\R\chi(x)\,dx_1
\le \left(\chi_1(\hat x_1)\int_\R\chi_2(\hat x_2)\,
dx_1\dots\int_\R\chi_n(\hat
x_n)\,dx_1\right)^{1/(n-1)}.$$
So, integrating in~$dx_2$ and using~\eqref{oh}
and~\eqref{oh2} once more,
\begin{eqnarray*}&&
\int_{\R^2} \chi(x)\,d(x_1,x_2)
\\ &\le& \left(\int_\R \chi_1(\hat x_1)\,dx_2 \int_\R\chi_2(\hat x_2)
dx_1\right.
\\ && \left. \qquad \int_{\R^2}\chi_3(\hat x_3)\,d(x_1,x_2)
\dots\int_{\R^2} \chi_n(\hat
x_n)\,d(x_1,x_2)\right)^{1/(n-1)}\\ &=&
\left(\int_\R \chi_1(\hat x_1)\,d\hat x_{1;2} \int_\R\chi_2(\hat x_2)
d\hat x_{2;2} \right.\\
&&\left.\qquad
\int_{\R^2}\chi_3(\hat x_3)\,d(x_1,x_2)
\dots\int_{\R^2} \chi_n(\hat
x_n)\,d(x_1,x_2)\right)^{1/(n-1)}
,
\end{eqnarray*}
where we denoted by~$d(x_1,x_2)$ the volume element in~$\R^2$.
By iterating this argument, for any~$k\le n$, we conclude that
\begin{eqnarray*}&&
\int_{\R^k} \chi(x)\,d(x_1,\dots,x_k)
\\ &\le&
\left(\int_{\R^{k-1}} \chi_1(\hat x_1)\,d\hat x_{1;k}\dots
\int_{\R^{k-1}} \chi_k(\hat x_k)
d\hat x_{k;k} \right.\\
&&\left.\qquad
\int_{\R^k}\chi_{k+1}(\hat x_{k+1})\,d(x_1,\dots,x_k)
\dots\int_{\R^k} \chi_n(\hat
x_n)\,d(x_1,\dots,x_k)\right)^{1/(n-1)}
,
\end{eqnarray*}
and, finally,
\begin{eqnarray*}
\int_{\R^n} \chi(x)\,dx
&\le&
\left(\int_{\R^{n-1}} \chi_1(\hat x_1)\,d\hat x_{1;n}\dots
\int_{\R^{n-1}} \chi_n(\hat x_n)
d\hat x_{n;n}
\right)^{1/(n-1)}
.
\end{eqnarray*}
This and~\eqref{oh3}
imply the claim in~(i).
Then,~(ii) easily follows from~(i).
\end{proof}

As a curiosity, we remark that the estimates in Lemma~\ref{proj}
are optimal, as the example of the cube shows, and that they
may be seen as suitably refined versions of
the classical isoperimetric and isodiametric
inequalities.

The main estimate needed for the proof of Theorem~\ref{GMT}
is the following:

\begin{lemma}\label{gmt}
Let $s\in (0,1)$.
Let $A$ and $B$ be disjoint measurable subsets of $\R^n$,
with $|A|=1$. Let $D:=\CC(A\cup B)$. Then, there exists
$\delta \in(0,1/10)$ depending on $n$ and $s$ such that the following
holds:
if $|B|\le \delta$, then
$$ L(A,D)\ge \left\{
\begin{matrix}
\delta &\qquad{\mbox{ if $s\in(0,1/2)$,}}  \\
\delta \log ({1}/{|B|}) &\qquad{\mbox{ if $s=1/2$,}}\\
\delta \,|B|^{1-2s}&\qquad{\mbox{ if $s\in(1/2,1)$.}}
\end{matrix}
\right.$$
\end{lemma}

\begin{proof}

The main step of
the proof consists in the following estimate: there
exists $\tilde c \in(0,1)$, suitably small,
depending on $n$ and $s$, such that,
for any \begin{equation}\label{eta}
r \in \left[\tilde C |B|,\, \tilde c \right],
\end{equation}
with $\tilde C:=1/\tilde c$, we have
\begin{equation}\label{corner}
\int_A \left[ \int_{D \cap (B_{\tilde Cr}(x)\setminus B_{\tilde c r}(x))}
\frac{dy}{|x-y|^{n+2s}}\right]\,dx\ge \tilde c\,{r^{1-2s}}.
\end{equation}

In order to prove \eqref{corner},
we divide $\R^n$ into a collection $K$ of
nonoverlapping cubes $Q$
of size $r$. We define
\begin{eqnarray*}
&& K_B:=\left\{ Q\in K {\mbox{ s.t. }} \frac{|Q\cap B|}{|Q|} \ge \frac 13
\right\},\\
&& K_D:=\left\{ Q\in K\setminus K_B {\mbox{ s.t. }} \frac{|Q\cap
D|}{|Q|}
\ge \frac 13 \right\},\\
&& K_A:= K\setminus (K_B\cup K_D).
\end{eqnarray*}

We also set
$$ Q_B:=\bigcup_{Q\in K_B} Q,\qquad\quad
Q_D:=\bigcup_{Q\in K_D} Q,\quad\qquad
Q_A:=\bigcup_{Q\in K_A} Q.$$
We observe that
\begin{equation}\label{BqB}
\delta \ge
|B|\ge \sum_{Q\in K_B} |Q\cap B|\ge \frac 1 3 \sum_{Q\in K_B} |Q|= \frac 13 \,{|Q_B|}.
\end{equation}
Moreover,
\begin{equation}\label{55}
\begin{split}
{\mbox{if $Q\in K_A$, }}& |Q\cap A|=|Q|-|Q\cap B|-|Q\cap D|\\
&\qquad \ge |Q|-\frac 1 3 |Q|-\frac 13|Q| = \frac 1 3 |Q|.\end{split}\end{equation}
In particular,
\begin{equation}\label{6712}
|Q_A| =\sum_{Q\in K_A}|Q|\le 3 \sum_{Q\in K_A} |Q\cap A|\le
3|A|=3.
\end{equation}
We also point out that if $x\in Q\subset K_D$ and $\tilde C> \sqrt n$,
then $Q\subseteq B_{\tilde C r}(x)$ and so
\begin{equation}\label{12ca}
\begin{split}
& \int_{Q\cap A} \left[ \int_{D \cap (B_{\tilde Cr}(x)\setminus
B_{\tilde cr}(x))}
\frac{dy}{|x-y|^{n+2s}}\right]\,dx
\\ &\quad\ge
\int_{Q\cap A} \left[ \int_{(Q\cap D)\setminus B_{\tilde cr}(x)}
\frac{dy}{|x-y|^{n+2s}}\right]\,dx \\&\quad
\ge\int_{Q\cap A}\; \frac{|(Q\cap D)\setminus B_{\tilde cr}(x)|}{(\sqrt n r)^{n+2s}} \,dx
\\&\quad
\ge\frac{|Q\cap A|\; \Big(r^n/3-
\tilde c^n\,|B_1|\, r^n\Big)}{(\sqrt n r)^{n+2s}} \ge c r^{-2s} |Q\cap A|,
\end{split}
\end{equation}
provided that $\tilde c$ is sufficiently small.
Now, two cases may occur. Either
\begin{equation}\label{1ca}
|Q_D\cap A|\ge r
\end{equation}
or not.
If \eqref{1ca} holds, then we exploit \eqref{12ca} to obtain
\begin{eqnarray*}
&& \int_A \left[ \int_{D \cap (B_{\tilde C r}(x)\setminus B_{\tilde c r}(x))}
\frac{dy}{|x-y|^{n+2s}}\right]\,dx
\\ &&\qquad\ge \sum_{Q\in K_D}\int_{Q\cap A} \left[ \int_{D \cap
(B_{\tilde C r}(x)\setminus
B_{\tilde c r}(x))}
\frac{dy}{|x-y|^{n+2s}}\right]\,dx
\\ &&\qquad\ge
\sum_{Q\in K_D}cr^{-2s}\,|Q\cap A| \ge cr^{-2s}|Q_D \cap A| \ge cr^{1-2s}.
\end{eqnarray*}
That is, if \eqref{1ca} holds true, then \eqref{corner}
is proved, up to renaming the constants.

Therefore, we can focus on the
case in which \eqref{1ca}
does not hold and suppose from now on that
\begin{equation}\label{2ca}
|Q_D\cap A|< r.
\end{equation}
Hence, recalling \eqref{eta} and~\eqref{BqB},
and the fact that~$\delta<1/10$, we conclude that
\begin{equation}\label{p1}
\begin{split}& |Q_A|\ge |Q_A\cap A|=|A|-|Q_B\cap A|-|Q_D\cap A|\\
&\quad\ge 1-|Q_B|- r \ge \frac12.\end{split}\end{equation}
{F}rom \eqref{p1} and Lemma \ref{proj}(ii), we have that, up to rotation,
\begin{equation}\label{78p}
|\Pi_{e_n} (Q_A)|\ge c_0.
\end{equation}
Thus, we organize the cubes of $K$ into subfamily
of columns in direction $e_n$:
more explicitly, the column containing
a cube $Q\in K$ is given by the union of all the cubes
of the form $Q+j \eta e_n$, for any $j\in \Z$.

We define $C_A$ to be the union of all the columns
that have a cube belonging to $K_A$,
and $C_B$ to be the union of all the columns in $C_A$
that have at least one cube belonging
to $K_B$\footnote{Of course, $C_A$ and $C_B$
may well have some common columns in the
intersection.}. We also let $M_A$
and $M_B$ to be the cardinality of the columns belonging to
$C_A$ and to $C_B$, respectively. We remark that
\begin{equation}\label{af}{\mbox{
the number of cubes
in~$K_A$ is finite,}}\end{equation} due to~\eqref{6712},
and so $M_A$ 
is finite\footnote{Similarly, exploiting \eqref{BqB},
one can see that~$M_B$ is finite -- however,
a more precise estimate on~$M_B$ will be
given in the forthcoming~\eqref{mb}.}.

Notice that $Q_A\subseteq C_A$, therefore, by \eqref{78p},
\begin{equation}\label{76}
c_0 \le |\Pi_{e_n} (C_A)|=r^{n-1} M_A.
\end{equation}

On the other hand, if $C_o$ is a column belonging to $C_B$, then
it contains one cube $Q^{(1)}$
belonging to $K_B$, and therefore
$$|C_o \cap B| \ge |Q^{(1)} \cap B|\ge r^n/3.$$
Consequently,
$$ |B| \ge |C_B \cap B|\ge M_B r^n /3.$$
Accordingly, recalling \eqref{eta} and \eqref{76},
\begin{equation}\label{mb}
 M_B \le\frac{3\,|B|}{ r^n}\le
\frac{3 \tilde c r}{r^{n}}\le \frac{c_0}{2r^{n-1}} \le \frac{M_A}2,
\end{equation}
if $\tilde c \le c_0/6$.
As a consequence of this, using \eqref{76} once more, we conclude that
\begin{equation}\label{CACB}
\begin{split}&{\mbox{the number of columns in $C_A\setminus C_B$}}\\
&\quad{\mbox{
is at least }} \,M_A-M_B\ge \frac{M_A}2\ge \frac{c_0}{2r^{n-1}}.
\end{split}\end{equation}
Now, let $C^\star$ be a column in $C_A\setminus C_B$. Then, recalling
\eqref{af},
we see that $C^\star$ must contain only a finite number of cubes
belonging to $K_A$, so we may define the cube $Q^\star_o$
as the cube of $C^\star$ belonging to $C_A$ with the
highest possible $e_n$-coordinate.

We consider the cube
$Q^\star_1:=Q^\star_o+2re_n$. By construction $Q^\star_1 \in K_D$. Notice that if $x\in Q^\star_o$ and $y\in Q^\star_1$ then
$|x-y|\ge r$, and
$|x-y|\le (2 + \sqrt n ) r \le \tilde C r$, provided that $\tilde C$ is 
sufficiently large. Therefore, if~$x\in Q^\star_o$, 
then~$Q^\star_1\subseteq
B_{\tilde C r}(x)
\setminus
B_{\tilde c r}(x)$, and
\begin{eqnarray*}&&
\int_{A\cap Q^\star_o} \int_{D \cap
(B_{\tilde C r}(x)
\setminus
B_{\tilde c r}(x))}
\frac{dx\,dy}{|x-y|^{n+2s}}\\ &&\quad\ge
\int_{A\cap Q^\star_o} \int_{D \cap Q^\star_{1}
}\frac{dx\,dy}{|x-y|^{n+2s}}
\\&&\quad\ge c r^{-(n+2s)}\,|A\cap Q^\star_o|\,
|D \cap Q^\star_{1}| \\&&\quad
\ge
c r^{n-2s},\end{eqnarray*}
for a suitable~$c\in(0,1)$ (independent of~$\tilde c$ and~$\tilde C$,
and possibly different line after line).
As a consequence,
$$ \int_{C^\star\cap A} \int_{ D
\cap (B_{\tilde C r}(x) \setminus B_{\tilde c r}(x)) }
\frac{dx\,dy}{|x-y|^{n+2s}}
\ge c r^{n-2s}.$$
Since this is valid for any
column $C^\star$ in $C_A\setminus C_B$, in the light of \eqref{CACB}
we obtain that
$$ \int_{A}
\int_{D \cap (B_{\tilde C r}(x) \setminus B_{\tilde c r}(x))}
\frac{dx\,dy}{|x-y|^{n+2s}}
 \ge c r^{1-n}r^{n-2s}=cr^{1-2s}
.$$
This completes the proof of~\eqref{corner}, by choosing $\tilde c
\in(0,1)$ small enough.

Since~$\tilde c$ is now fixed once and for all,
we can suppose that~$\delta$ in the statement of Lemma~\ref{gmt} is
smaller than~$\tilde c^3$, hence~$\tilde c^3\ge |B|$. So,
let $k_0$ be the largest integer so that $\tilde c^{2k_0+1} \ge |B|$.
{F}rom~\eqref{eta} and~\eqref{corner}, we deduce that,
for any $1 \le k \le k_0$
$$   \tilde C |B| \le \tilde r_k:=\tilde c ^{2k} \le \tilde c ,$$
we have that
$$
\int_A \left[ \int_{D \cap (B_{\tilde r_k/\tilde c}(x)\setminus
B_{\tilde c r_k}(x))}
\frac{dy}{|x-y|^{n+2s}}\right]\,dx\ge \tilde c^{1+2k(1-2s)}.$$
Consequently,
\begin{eqnarray*}
&& \int_A \left[ \int_{D}
\frac{dy}{|x-y|^{n+2s}}\right]\,dx
\\&\ge& \sum_{k=1}^{k_0}
\int_A \left[ \int_{D \cap (B_{\tilde r_k/\tilde c}(x)\setminus
B_{\tilde c r_k}(x))}
\frac{dy}{|x-y|^{n+2s}}\right]\,dx\\
&\ge& \sum_{k=1}^{k_0} \tilde c^{1+2k(1-2s)}\\ &\ge&
\left\{
\begin{matrix}
\tilde c^{2 (1-2s)+1} 
& {\mbox{ if $s\in(0,\,1/2)$,}}\\
k_0\tilde c & {\mbox{ if $s=1/2$,}}\\
\tilde c^{2k_0(1-2s)+1}
& {\mbox{ if $s\in(1/2,\,1)$,}}
\end{matrix}
\right.
\end{eqnarray*}
This implies the desired result, by taking~$\delta$
appropriately small with respect to~$\tilde c$.
\end{proof}

We now complete the proof of Theorem~\ref{GMT}
by the following argument. If~$|B|>c |A|$, then~$|A\cup B|\le
|A|+|B|<(1+c^{-1})|B|$. Hence, we make
use of the Sobolev-type inequality in~\eqref{CA3}
and we possibly rename~$c$, to
conclude that
\begin{equation}\label{HERE}\begin{split}
& L(A,D) =
\int_{A} \int_{\CC(A\cup B)}
\frac{dx\,dy}{|x-y|^{n+2s}}\\ &\qquad\ge
c\,|A||A\cup B|^{-2s/n} \ge c\,|A||B|^{-2s/n}
\end{split}\end{equation}
which proves~\eqref{gmt2}.

If, on the other hand, if~$|B|\le c |A|\ne0$,
we define
$$ \tilde A:= \frac{1}{|A|^{1/n}} A, \;\tilde B:= \frac{1}{|A|^{1/n}} B
{\mbox{
and }}\;
\tilde D:=\frac{1}{|A|^{1/n}}D=\CC(\tilde A\cup\tilde B).$$
We observe that~$|\tilde A|=1$ and~$|\tilde B|=|B|/|A|\le c$,
so we can apply Lemma~\ref{gmt}
and obtain that
\begin{eqnarray*}
L(A,D)&=& |A|^{(n-2s)/n} L(\tilde A,\tilde D)\\
&\ge& \left\{
\begin{matrix}
c\,|A|^{(n-2s)/n}
&\qquad{\mbox{ if $s\in(0,1/2)$,}}\\
c\,|A|^{(n-1)/n}
\log ({1}/{|\tilde B|}) &\qquad{\mbox{ if $s=1/2$,}}\\
c\,|A|^{(n-2s)/n}
\,|\tilde B|^{1-2s}&\qquad{\mbox{ if $s\in(1/2,1)$.}}\\
\end{matrix}
\right.\\ &=&
\left\{
\begin{matrix}
c\,|A|^{(n-2s)/n}
&\qquad{\mbox{ if $s\in(0,1/2)$,}}\\
c\,|A|^{(n-1)/n} \log ({|A|}/{|B|}) &\qquad{\mbox{ if $s=1/2$,}}\\
c\,|A|^{(n-2s)/n} \,
(|B|/|A|)^{1-2s}&\qquad{\mbox{ if $s\in(1/2,1)$.}}\\
\end{matrix}
\right.
\end{eqnarray*}
This proves~\eqref{gmt2} and completes the proof of
Theorem~\ref{GMT}.

At the end of this section we prove another localized version
of~\eqref{gmt2} in Theorem ~\ref{GMT} which could be useful in other situations,
such as in \cite{SaV}.

\begin{prop}\label{GMTloc}
Let $ s\in [1/2,1)$. Let $A$, $D$ be disjoint subsets of a cube $Q
\subset \R^n$ with $$\min\{|A|,|D|\} \ge \sigma |Q|,$$ for some $\sigma
>0$. Let $B=Q\setminus(A\cup D)$. Then,
$$ L(A,D)\ge \left\{
\begin{matrix}
\delta |Q|^\frac{n-1}{n}\,\log (|Q|/{|B|}) &\qquad{\mbox{ if $s=1/2$,}}\\
\ \\
\delta |Q|^\frac{n-2s}{n}\,(|Q|/|B|)^{2s-1}&\qquad{\mbox{ if
$s\in(1/2,1)$.}}\\
\end{matrix}
\right.$$
for some $\delta >0$ depending on $\sigma$, $n$ and $s$.
\end{prop}

\begin{proof}
By rescaling, we can assume that $Q$ is the unit cube.

The proof follows the lines of Lemma \ref{gmt}: we show that for each $$r \in [\tilde C |B|,\tilde c],$$
we satisfy \eqref{corner} with~$\tilde C:=1/\tilde c$, for some small $\tilde c$ depending on $n$, $s$ and $\sigma$.

We divide the unit cube into cubes of size $r$ that we partition into the three sets $K_A$, $K_B$ and $K_D$. The difference now is that $D$ is defined only
on the unit cube and therefore the existence of $Q^*_1 \in K_D$ at the end of the proof of Lemma \ref{gmt} requires a more careful argument.

As before, we only need to deal with the case (see \eqref{p1})

$$|Q_A| \ge \frac \sigma 2, \quad |Q_D| \ge \frac \sigma 2.$$ Denote $$\alpha:=|Q_A|, \quad \frac \sigma 2 \le \alpha \le 1- \frac \sigma 2.$$

{F}rom Lemma \ref{proj}(ii), we have that, up to rotation,
\begin{equation}\label{78p1}
|\Pi_{e_n} (Q_A)|\ge \alpha^\frac{n-1}{n}.
\end{equation}

We define $C_A$ to be the union of all the columns
that have a cube belonging to $K_A$,
$C_B$ to be the columns in $C_A$
that have at least one cube belonging
to $K_B$ and $C_D$ to be the columns in $C_A \setminus C_B$ that have at least one cube belonging to $K_D$. We also let $M_A$
$M_B$, and $M_D$ to be the cardinality of the columns belonging to
$C_A$, $C_B$ and to $C_D$, respectively.

By \eqref{78p1},
\begin{equation*}
\alpha^\frac{n-1}{n} \le |\Pi_{e_n} (C_A)|=r^{n-1} M_A.
\end{equation*}
and we also have (see \eqref{mb}) $$M_B \le 3 \tilde c r^{1-n}.$$
On the other hand, the cardinality $m_a$ of the columns in $C_A \setminus (C_B \cup C_D)$ (that contain only cubes from $K_A$) satisfies
$$m_a r^{n-1} \le |Q_A| = \alpha.$$
Thus, if $\tilde c$ is sufficiently small, $$M_D=M_A-m_a - M_B \ge r^{1-n} (\alpha^ \frac{n-1}{n} - \alpha -\tilde c /3) \ge c r^{1-n}.$$

Let $C^*$ be a column belonging to $C_D$. Then $$C^* \subset K_A \cup K_D, \quad C^*\cap K_A \ne \emptyset, \quad C^*\cap K_D\ne \emptyset,$$
and we can easily conclude that on this column there must exist two cubes $Q^*_0 \subset K_A$, $Q^*_1 \subset K_D$ at distance either $2r$ or $3r$ from each other. Then $$L(A\cap Q_0^\star, D \cap Q_1^\star) \ge cr^{n-2s},$$ and the proof continues as before.
\end{proof}

\appendix

\section{A Sobolev-type inequality for sets}\label{APP}

For completeness, we give here
an elementary proof of
the Sobolev-type inequality exploited 
in our paper,
see \eqref{HERE}.
For related and more general results see \cite{BBrMi, 
BrMi} and 
references
therein.

\begin{lemma}\label{5yhh}
Fix $x\in\R^n$. Let $E\subset\R^n$ be a measurable set with finite
measure.
Then,   
$$ \int_{\CC E}\frac{dy}{|x-y|^{n+2s}}\ge c(n,s)\,|E|^{-2s/n},$$
for a suitable constant $c(n,s)>0$.
\end{lemma}

\begin{proof} Let
$$ \rho:=\left(\frac{|E|}{\omega_n}\right)^{1/n}.$$
Then,
\begin{eqnarray*}&& |\CC E\cap B_\rho(x)|=|B_\rho(x)|- |E\cap B_\rho(x)|
\\&&\qquad=
|E|-|E\cap B_\rho(x)|=|E\cap \CC B_\rho(x)|.\end{eqnarray*}
Therefore,
\begin{eqnarray*}
&& \int_{\CC E}\frac{dy}{|x-y|^{n+2s}}=
\int_{\CC E\cap B_\rho(x)}\frac{dy}{|x-y|^{n+2s}}+
\int_{\CC E\cap \CC B_\rho(x)}\frac{dy}{|x-y|^{n+2s}}
\\ &&\qquad \ge\int_{\CC E\cap B_\rho(x)}\frac{dy}{\rho^{n+2s}}+
\int_{\CC E\cap \CC B_\rho(x)}\frac{dy}{|x-y|^{n+2s}}
\\ &&\qquad
=\frac{|\CC E\cap B_\rho (x)|}{\rho^{n+2s}}+
\int_{\CC E\cap \CC B_\rho(x)}\frac{dy}{|x-y|^{n+2s}}
\\ &&\qquad=
\frac{|E\cap \CC B_\rho (x)|}{\rho^{n+2s}}+
\int_{\CC E\cap \CC B_\rho(x)}\frac{dy}{|x-y|^{n+2s}}
\\ &&\qquad\ge
\int_{E\cap \CC B_\rho(x)}\frac{dy}{|x-y|^{n+2s}}
+
\int_{\CC E\cap \CC B_\rho(x)}\frac{dy}{|x-y|^{n+2s}}
\\ &&\qquad
=\int_{\CC B_\rho(x)} \frac{dy}{|x-y|^{n+2s}}.
\end{eqnarray*}
By using polar coordinate centered at $x$,
the desired result easily follows.
\end{proof}

By integrating in $x$ the estimate of Lemma \ref{5yhh},
we obtain:

\begin{corollary}
Let $E$, $F\subset\R^n$ be measurable sets with finite measure.
Then,
\begin{equation}\label{CA3}\int_F \int_{\CC
E}\frac{dx\,dy}{|x-y|^{n+2s}}\ge
c(n,s)\,|F|\,|E|^{-2s/n},\end{equation}
for a suitable constant $c(n,s)>0$.

In particular,
\begin{equation}\label{CA3bis}\int_E \int_{\CC
E}\frac{dx\,dy}{|x-y|^{n+2s}}\ge
c(n,s)\,|E|^{(n-2s)/n}\end{equation}
for any
measurable set~$E$ with finite measure.
\end{corollary}

\vfill

\vspace{1cm}

{{\sc Ovidiu Savin}

Mathematics Department, Columbia University,

2990 Broadway, New York , NY 10027, USA.

Email: {\tt savin@math.columbia.edu}
}

\vspace{1cm}

{{\sc Enrico Valdinoci}

Dipartimento di Matematica, Universit\`a di Roma Tor Vergata,

Via della Ricerca Scientifica 1, 00133 Roma, Italy.

Email: {\tt enrico@mat.uniroma3.it}
}

\end{document}